\documentclass[np,hyp,12pt]{nyj}
\begin{document}
\thispagestyle{empty}
\title{Electronic Mathematics Journals}
\author{Mark Steinberger}
\address{Department of Mathematics and Statistics\\ The University at
Albany, State University of New York\\ Albany, NY 12222}
\email{mark@csc.albany.edu}

\begin{abstract} In the Forum section of the November, 1993 {\em
Notices of the American Mathematical Society\/},
John Franks discussed the electronic journal of the future. 
Since then, the New York Journal of Mathematics, the first
electronic general mathematics journal, has begun publication. 
In this article, we 
explore the issues of electronic journal publishing in the context of
this new project. We also discuss future developments.
\end{abstract}
\maketitle
\tableofcontents

John Franks, in the November, 1993 Forum section 
of the {\em Notices of the American Mathematical Society\/}
\cite{F}, stressed
the advantages of nonprofit electronic
journals as a vehicle for mathematical research. We agree.

\vspace{10pt}

The primary advantage for most mathematicians is ease of access.
Without leaving one's desk, one may browse the articles and print out
any articles deserving of more detailed consideration.

 From the standpoint of the author, the delays in publication can be
confined to the refereeing process, with approved articles appearing
shortly after the peer reviewing process is complete. 

 From an institutional point of view, there is a big
financial advantage in promoting nonprofit electronic publications.
The cost of producing, distributing, and archiving
an electronic journal is smaller than
that for a print journal (cf.~\cite{Od,Q1}. The amount of savings
available through electronic publishing has been debated, but all
commentators seem to agree the costs are significantly lower.)
In a time of declining resources and
escalating prices for print journals, this is an important advantage.

There are also additional capabilities available with electronic
media. (See Section~\ref{enhancsec} for details, including information
regarding the implementation of these capabilities in the New
York Journal of Mathematics.) 
For instance,
abstracts may be distributed over mailing lists, notifying the
reader of the availability of the articles. Electronic links may be
made to past and future papers, reviews, comments, and elucidations of
the work in question. And keyword searches may be made to identify
articles of interest to researchers. 
These capabilities vastly increase the ability of the author
and the journal to provide information to the reader. 

\vspace{10pt}

At the New York Journal of Mathematics, we are
in process of working out the practical side of the issues surrounding
electronic publication.
We understand that the formats and functions of electronic journals
will evolve over time. We feel it is essential that this evolution
take place in a manner that will maximize quality and work to the
maximum benefit of the mathematical community.

\section{Editorial Issues}
We feel the quality of the mathematics in our
articles is of
utmost importance. We have implemented the traditional peer review
process in its full rigor. Papers are blind-refereed for quality and
correctness, as is done for high quality print journals.

 From the standpoint of the issues surrounding
electronic publication, we have taken the following considerations to
be fundamental.

\vspace{4pt}

\begin{itemize}
\item The appearance of every paper printed out from our journal should be
uniform and appealing.

\vspace{10pt}

 Thus, we decided not to distribute ascii source files, but to
distribute our papers in graphical formats only. The papers
   are typeset in a traditional format, in accordance with the
   journal's style sheet, and with the logo, the statement of
   copyright, ISSN, and pagination given on the first page.

\vspace{10pt}

\item The standards for succinctness in writing should be consonant with
those applied by the other journals in the field.

\vspace{10pt}

 Some proponents of electronic publication have urged changes in style,
citing the low price of disk space as a rationale for publishing
articles more loquacious than those
commonly acceptable in a print medium. We
decided to eschew this route, on the grounds that the perceived
quality of our publications would be reduced. We feel it is important
to follow the standards of consensus in the field. If these standards
change in the future, we will change with them.

\vspace{10pt}

Authors who wish to expand on the material in greater detail than the
editor feels is appropriate for a journal article are welcome to
submit additional material for inclusion in hypertext comment files
(see below).

\vspace{10pt}

\item The written record must be maintained intact in perpetuity.

\vspace{10pt}

The University at Albany,
State University of New York, has endorsed this commitment, agreeing
to insure the integrity of the journal's
archive in perpetuity.

Our articles are fixed at time of publication, with their pages
numbered consecutively throughout each volume, in the traditional
manner.

\end{itemize}

\vspace{6pt}

These decisions have affected some of our further options, while
leaving others open.

Most prominently, the decision not to distribute ascii source files
has made distribution by email impractical. Our papers are available
through internet tools such as ftp, gopher, and the World Wide Web (WWW).

\section{Electronic Features and Enhancements}\label{enhancsec}
\subsection*{Client and Server Technologies}
As is indicated by the variety of the internet tools
in use, internet technology is
a moving target, and electronic journals must move with it. 
In particular, the journals must provide service to accommodate the
technologies used by the various users.

Ftp and gopher were developed earlier than WWW, and
software for accessing materials served by ftp and gopher is still
more common (but not by much) than software for accessing materials
over the World Wide Web. However, the capabilities of the web are
greater, as its servers and clients can communicate via 
{\em hypertext\/},
a system in which access to files is controlled by electronic links
embedded in text files.

This provides a distinct advantage over the pre-existing technologies,
which basically provide access to directories of files. With
hypertext, one may offer the reader a direct link to another resource
anywhere on the internet, from any point in the file currently being
accessed. By activating the link, the reader is given immediate access
to the resource in question, and may return to 
precisely the same location in preceding document
with a keystroke or the click of a button. This makes
cross-referencing very flexible and effective, increasing the
efficiency of searching out information considerably.

It is also useful to be able to decouple the presentation of
information from a hierarchical, tree-like structure, and allow
authors and publishers to connect material as dictated by its internal
logic. 

\subsection*{Hypertext features useful in electronic journals}
Hypertext is particularly well adapted to cross-references between
papers, or between papers and reviews and other discussions. 
And those journals offered via the World Wide Web can take advantage
of this.

At the simplest level of implementation,
each paper can have a bibliography file
containing links to any item in the paper's bibliography
available on the net. And comment files are especially suited to
hypertext format.

Comment files are an important innovation available
to electronic journals. Maintained by the editor, they may include
links to reviews, to subsequent articles in which the results are
extended or applied, to errata, to elucidations by the author of material in
the paper, etc. Using hypertext, such files can be structured as
narratives together with links to the resources they reference.

Note that the involvement of the editor in screening the material to
be placed in the comment file is a useful filter against
clutter in the literature. Indeed, while the New York Journal has
advertised its intention to provide such files, we have not yet
received submissions for them.

An additional level of connectedness is given by the advent of
reviewing journals on the World Wide Web. Both {\em Math Reviews\/}
and {\em Zentralblatt\/} are soon to be offered on the web. In the
case of {\em Math Reviews\/}, the reviews will include direct links to
papers offered by electronic journals, as well as numerous internal
links to other reviews. Hopefully, it will also be possible to make
links from comment files in electronic journals to reviews in the
reviewing journals.

This opens up a very interesting prospect, whereby it is possible to
follow an idea through several different papers and reviews via
seemless, near-instantaneous linkages, without ever getting up from
one's desk. This is likely to be a great boon for scholarship.

\subsection*{Embedding hypertext links directly in journal articles}
Just as hypertext itself is more flexible and useful than offering
simple directories of files, it is much more useful to be able to
embed hypertext links directly in the graphical formats of the papers
themselves, rather than restricting them to auxiliary files associated
to the papers.

The technology to do
this in \TeX{} documents is now becoming available through the
Hyper\TeX{} project, an offshoot of Paul Ginsparg's e-print server
project at Los Alamos National Labs. 

Hyper\TeX{} provides a method of embedding hypertext links directly in dvi
files, postscript files, or Acrobat's pdf files. These links can be
followed if you make use of a viewer programmed to recognize and
follow them. Such viewers are either available or about to become available
on all major computing platforms. The interested reader may find out
more about this project from the following locations on the World Wide
Web:
\begin{center}
\leavevmode\url{http://nyjm.albany.edu:8000/hyper.html}\\
\leavevmode\url{http://math.albany.edu:8800/hm/ht/}
\end{center}

Hyper\TeX{} links have another benefit: internal cross-references.
Links for the internal cross-references in a paper allow the
reader to quickly flip to the statement of a theorem when it is
invoked, supplying the information needed to understand an argument.
Then, with a keystroke, the reader may flip back to the argument and
continue reading. The process is much faster than flipping pages in a
paper copy of the article, giving a real advantage to reading with a
\TeX{} viewer.

And a table of contents for the article, including links to the various
sections, is quite useful for browsing.

The New York Journal of Mathematics offers Hyper\TeX{} dvi files for all
papers published since March, 1995.

\subsection*{Keyword searching} Through a technique called WAIS
indexing, it is possible to build databases of text from internet
sites and to run keyword searches of those databases. This can be set
up to produce electronic links to those documents that match the
keywords. 

It is also possible to index the full \TeX{} source of articles, and
then pass links either to graphical formats or to a ``home page'' for
the paper when a match occurs. The home page can contain links to
various graphical formats.
Such a schema has been implemented for the New York Journal of Mathematics at
\begin{center}
\leavevmode\url{http://nyjm.albany.edu:8000/SF/nyjmsearch.html}
\end{center}
In the New York Journal, a paper's home page is a WWW page 
containing the paper's
abstract, keywords, and subject classification, along with various links.

While WAIS technology is useful for searching a particular 
internet site, it
can be even more useful to search an index of many different sites.
Such an index exists in Australia: Jim Richardson's MathSearch
index, at 
\begin{center}
\leavevmode\url{http://ms.maths.usyd.edu.au:8000/MathSearch.html}
\end{center}
It indexes all the major mathematical sites on the web, worldwide.
A keyword search on this database can turn up links to resources
at sites one didn't even know existed.

\subsection*{Direct communication to the readers} Electronic
communication permits direct distribution of abstracts through
electronic mail to interested readers. This notifies the reader of the
existence of a paper and permits him or her to fetch the full paper if
interested. 

The New York Journal of
Mathematics maintains four listserv lists for this purpose, running on
listserv@albany.edu. 
One list, nyjmth-a,
distributes abstracts for all papers. The other three
are specialty lists, distributing
abstracts in algebra, analysis, and geometry/topology. Their list
names are nyjm-alg, nyjm-an, and nyjm-top, respectively.

In future, it should be possible to design systems tailored more
precisely to the specific interests of the readers. For instance,
readers could be
sent abstracts of all papers containing reader-specified subject
classification numbers in the author's list of primary and secondary
subject classifications.

\section{The New York Journal of Mathematics}
The New York Journal of Mathematics was launched with the support and
assistance of
the Office of
Information Systems and Technology of the University at Albany, State
University of New York.
We also have ongoing support from 
the Office of the Vice 
President for Research and the
Department of Mathematics and Statistics.

The Office of Information Systems and Technology includes both the
University Libraries and Computing and Network Services. These two
units have collaborated in sponsoring the
University's Electronic Library Initiative. The journal itself is part
of this initiative. 

The journal is integrated with the library system, and
professional librarians have been active in every phase of the
journal's development, including the provisions for the preservation
of the data and the maintenance of the integrity of the archive.
Indeed, the New York Journal complies with the
recommendations given in \cite{Q2}.

The New York Journal of Mathematics is available by WWW, gopher and
ftp. WWW access is given by
\begin{center} \leavevmode\url{http://nyjm.albany.edu:8000/nyjm.html}\end{center}
while gopher access is given by the command
\begin{center} gopher nyjm.albany.edu 1070\end{center}
The same archive is accessible by anonymous ftp on nyjm.albany.edu in
the directory /pub/nyjm.

Papers should be submitted by electronic mail, directly to the editor
whose field is closest to the work in question. Our editorial board is
listed at
\begin{center} \leavevmode\url{http://nyjm.albany.edu:8000/Edboard.html}\end{center}
Instructions for the preparation of articles may be found at
\begin{center} \leavevmode\url{http://nyjm.albany.edu:8000/Instr.html}\end{center}

Our papers are provided over the internet free of charge. Printed
copies may be obtained for a fee from the Department of Mathematics
and Statistics of the University at Albany, State University of New
York.


\begin{thebibliography}{9}
\bibitem{F} John Franks, {\em The impact of electronic
publication on scholarly journals\/}, Notices of the AMS {\bf 40}
(1993), 1200--1202. 

\bibitem{Od} Andrew Odlyzko {\em Tragic loss or good riddance? The
impending demise of traditional scholarly journals\/}, 
Notices of the AMS {\bf 42} (1995), 49--53.

\bibitem{Q1} Frank Quinn, {\em Roadkill on the electronic highway?
  The threat to the mathematical literature\/}, 
Notices of the AMS {\bf 42} (1995), 53--56.

\bibitem{Q2} Frank Quinn, {\em A role for libraries in electronic
publishing\/}, EJournal {\bf 4} No.~2 (June 1992).

\end{thebibliography}
\end{document}